\documentclass[a4paper]{article}
\usepackage{amsmath,amssymb,graphicx,mathrsfs,psfrag}

\DeclareMathOperator{\im}{Im}
\DeclareMathOperator{\cosech}{cosech}
\begin{document}
\title{A numerical method of Fourier transform based on hyperfunction theory}
\author{Hidenori Ogata\footnote{ 
Department of Computer and Network Engineering,  
Graduate School of Informatics and Engineering,  
The University of Electro-Communications,  
1-5-1 Chofu-ga-Oka, Chofu, 182-8585, Japan}}
\date{}
\maketitle
\begin{abstract}
In this paper, we propose a numerical method of Fourier transform based
 on hyperfunction theory. 
In the proposed method, we compute analytic functions 
called the defining functions, 
which give the desired Fourier transform as a hyperfunction, 
and then obtain the Fourier transform 
by the analytic continuation of the defining functions 
onto the real axis. 
Numerical examples show the efficiency of the proposed method compared to the previous methods.
\end{abstract}
\section{Introduction}
\label{sec:introduction}
Fourier transform
\begin{equation}
 \label{eq:Fourier-transform}
  \mathscr{F}[f](\xi) = 
  \int_{-\infty}^{+\infty}f(x)\mathrm{e}^{-2\pi\mathrm{i}\xi x}\mathrm{d}x
\end{equation}
is familiar in science and engineering. 
However, it is difficult to compute a Fourier transform by conventional
numerical integration formulas such as the DE rule \cite{TakahasiMori1974}, especially, 
if the integrand $f(x)$ decays slowly as $x\rightarrow\pm\infty$. 
We propose a numerical method for computing Fourier transforms
efficiently based on hyperfunction theory \cite{Sato1959}. 

Hyperfunction theory is a theory of generalized functions based on
complex analysis. 
Roughly speaking, a hyperfunction $f(x)$ is the difference of the boundary
values of an analytic function $F(z)$, that is, 
\begin{equation*}
 f(x) = F(x + \mathrm{i}0) - F(x - \mathrm{i}0), 
\end{equation*}
where $F(z)$ is called a defining function of the hyperfunction $f(x)$. 
For example, the Dirac delta function is defined by 
\begin{equation*}
 \delta(x) 
  = 
  - \frac{1}{2\pi\mathrm{i}}
  \left(
   \frac{1}{x + \mathrm{i}0} - \frac{1}{x - \mathrm{i}0}
  \right),
\end{equation*}
which coincides with an elementary definition of the delta function 
\begin{equation*}
 \delta(x) = 
  \lim_{\epsilon\rightarrow 0+}
  \frac{1}{\pi}\frac{\epsilon}{x^2+\epsilon^2}
\end{equation*}

In hyperfunction theory, the Fourier transform of a function $f(x)$ is
defined by 
\begin{align}
 \nonumber
 \mathscr{F}[f](\xi) = \: & 
 \int_{-\infty}^{0}f(x)\mathrm{e}^{-2\pi\mathrm{i}(\xi+\mathrm{i}0)x}\mathrm{d}x 
 + 
 \int_{0}^{\infty}f(x)\mathrm{e}^{-2\pi\mathrm{i}(\xi-\mathrm{i}0)x}\mathrm{d}x
 \\ 
 \label{eq:Fourier-transform2}
 = \: & 
 \mathfrak{F}_+(\xi+\mathrm{i}0) - \mathfrak{F}_-(\xi-\mathrm{i}0), 
\end{align}
where
\begin{align}
 \label{eq:defining-function1}
 \mathfrak{F}_+(\zeta) = \: & 
 \int_{-\infty}^{0}f(x)\mathrm{e}^{-2\pi\mathrm{i}\zeta x}\mathrm{d}x
 & ( \: \im\zeta > 0 \: ), 
 \\ 
 \label{eq:defining-function2}
 \mathfrak{F}_-(\zeta) = \: & 
 - \int_0^{\infty}f(x)\mathrm{e}^{-2\pi\mathrm{i}\zeta x}\mathrm{d}x
 & ( \: \im\zeta < 0 \: ).
\end{align}
It means that the Fourier transform $\mathscr{F}[f](\xi)$ is the
hyperfunction whose defining function is 
\begin{equation}
 \mathfrak{F}(\zeta) = 
  \begin{cases}
   \mathfrak{F}_+(\zeta) & ( \: \im\zeta > 0 \: ) \\ 
   \mathfrak{F}_-(\zeta) & ( \: \im\zeta < 0 \: ).
  \end{cases}
\end{equation}
In our method, we compute these defining functions $F_{\pm}(\zeta)$ in
$\mathbb{C}\setminus\mathbb{R}$ instead of the integral on the
right-hand side of (\ref{eq:Fourier-transform}), and we obtain the
Fourier transform by (\ref{eq:Fourier-transform2}) using the analytic
continuation of $\mathfrak{F}_{\pm}(\zeta)$ onto $\mathbb{R}$. 
It is easy to compute the integrals on the right-hand sides of 
(\ref{eq:defining-function1}) and (\ref{eq:defining-function2}) 
because their integrands decay exponentially as
$x\rightarrow\pm\infty$. 
The analytic continuation of $F_{\pm}(\zeta)$ is done using the
continued fraction expansions. 

Previous studies related to our paper are as follows. 
Toda and Ono proposed a method of computing Fourier integrals 
\begin{equation}
 \label{eq:Fourier-integral}
  I = \int_0^{\infty}f(x)\cos(2\pi\xi x+\theta_0)\mathrm{d}x, 
\end{equation}
where $\xi > 0$ and $\theta_0$ is a constant, 
by evaluating the limit
\begin{equation*}
 I = 
  \lim_{n\rightarrow\infty}\int_0^{\infty}
  f(x)\cos(2\pi\xi x+\theta_0)\exp(-2^{-n}x)\mathrm{d}x
\end{equation*}
using the DE rule and the Richardson extrapolation \cite{TodaOno1978}. 
Sugihara improved Toda and Ono's method 
and proposed a method by evaluating the limit
\begin{equation*}
 I =
  \lim_{n\rightarrow\infty}\int_0^{\infty}
  f(x)\cos(2\pi\xi x+\theta_0)\exp(-2^{-n}x^2)\mathrm{d}x, 
\end{equation*}
using the DE rule and the Richardson
extrapolation \cite{Sugihara1987}. 
Ooura and Mori proposed a DE-type numerical integration rule for oscillatory integrals 
with unique technique \cite{OouraMori1991}. 
They apply the variable transformation
\begin{equation*}
 x = \frac{1}{2\pi h}\varphi\left(u+\frac{h}{2}-\frac{h\theta_0}{2}\right), 
\end{equation*}
where 
\begin{equation*}
 \varphi(u) = \frac{u}{1-\exp(-2\pi \sinh u)}
\end{equation*}
and $h$ is a positive constant, to the integral
(\ref{eq:Fourier-integral}) and evaluate the transformed integral by
the trapezoidal rule with mesh $h$, that is, 
\begin{multline*}
 I \simeq 
 \frac{1}{2\xi}\sum_{k=-N_2}^{N_1}
 f\left(
 \frac{1}{2\xi h}
 \varphi\left(kh+\frac{h}{2}-\frac{h\theta_0}{\pi}\right)
 \right)
 \\ 
 \times
 \cos\left(
 \frac{\pi}{h}\varphi\left(kh+\frac{h}{2}-\frac{h\theta_0}{\pi}\right)
 + \theta_0
 \right)
 \varphi^{\prime}\left(kh+\frac{h}{2}-\frac{h\theta_0}{2}\right),
\end{multline*} 
where $N_1$ and $N_2$ are so small positive integers that the transformed integrand is very small at $k=N_1$ and $-N_2$. 
In the method, we can truncate the infinite sum of the trapezoidal rule  
with a small number of terms $N_1+N_2+1$ because the sampling points
rapidly approach the zeros of the integrand on the positive side of the $u$-axis, 
while the conventional DE rule are designed so that we can truncate the
trapezoidal rule with a small number of the sampling points by making
the transformed integrand decays double exponentially. 

As applications of hyperfunction theory to numerical analysis, 
Mori gave a theoretical analysis of numerical integration formulas based
on hyperfunction theory and showed that Gauss-type integration formulas
are obtained by approximating complex integrals which are defined as
hyperfunction integrals \cite{Mori1972}. 
The author proposed a numerical integration method based on
hyperfunction theory \cite{OgataHirayama2018}. 
In the paper, they obtain desired integrals by evaluating complex
integral which defines them as hyperfunction integrals, and they shows
the proposed method is efficient especially for integrals with strong
end-point singularities. 

The contents of the paper are as follows. 
In Section \ref{sec:hyperfunction}, we give a brief review of hyperfunction theory. 
In Section \ref{sec:numerical-Fourier-transform}, 
we propose a numerical method based on hyperfunction theory for computing Fourier transforms. 
In Section \ref{sec:example}, we give some numerical examples 
which show the effectiveness of the proposed method compared to the previous methods. 
In Section \ref{sec:concluding-remarks}, we give concluding remarks and refer to 
problems for future study. 
\section{Hyperfunction theory and numerical Fourier transform}
\label{sec:hyperfunction}
\subsection{Hyperfunctions}
We give a brief review of hyperfunction theory. 
For the detail of hyperfunctions, see \cite{Graf2010}.

Let $I$ be an open interval on the real axis $\mathbb{R}$ and $D$ be a
complex neighborhood of $I$, that is, a complex domain which includes
$I$ as a closed subset. Let $\mathscr{O}(D\setminus I)$ and
$\mathscr{O}(D)$ be the sets of all the holomorphic functions in
$D\setminus I$ and $D$ respectively. 
The set $\mathscr{O}(D\setminus I)$ forms a complex linear
space with the addition 
$F + G$ $( \: F, G \in \mathscr{O}(D\setminus I) \: )$ 
and the scalar product 
$cF$ $( \: F\in\mathscr{O}(D\setminus I), \: c\in\mathbb{C} \: )$  
respectively defined by 
\begin{equation*}
 (F+G)(z) := F(z) + G(z), \quad (cF)(z) := cF(z) \quad ( \: z\in
  D\setminus I \: ), 
\end{equation*} 
and the set $\mathscr{O}(D)$ can be regarded as a linear subspace of 
$\mathscr{O}(D)$. 

A hyperfunction $f=f(x)$ on the interval $I$ is an element 
(an equivalence class) of the quotient subspace%
\footnote{
As shown here, the behavior of the defining function $F(z)$ in the
vicinity of the interval $I$ is crucial for the hyperfunction $f=[F]$, and 
we do not need to cling to a particular choice of a complex neighbor
hood $D$ of $I$. 
Therefore, we should define the space of hyperfunctions on an interval
$I$ by the inductive limit 
\begin{equation*}
 \mathscr{B}(I) := \varinjlim_{D\supset I}
  \mathscr{O}(D\setminus I) / \mathscr{O}(D).
\end{equation*}
However, we do not go into the exact definition of $\mathscr{B}(I)$ any
more, and we do not need to do so \cite{Kaneko1988}.
}
\begin{equation*}
 \mathscr{B}(I):=\mathscr{O}(D\setminus I)\slash\mathscr{O}(D). 
\end{equation*}
If a hyperfunction $f$ is an equivalence class of a function 
$F\in\mathscr{O}(D\setminus I)$, we call $F$ a defining function of $f$,
which is denoted by 
\begin{equation*}
 f=[F] \quad \mbox{or} \quad f(x)=[F(z)]. 
\end{equation*}
We also denote a hyperfunction $f=[F]$ by 
\begin{equation*}
 f=[F_+, F_-] \quad \mbox{or} \quad f(x)=[F_+(z), F_-(z)],
\end{equation*}
where $F_{\pm}(z)$ is the restriction of $F(z)$ on 
\(
D_{\pm} := 
\left\{ \: z\in D \: | \: \pm\im z > 0 \: \right\}, 
\) 
that is, 
\begin{equation*}
 F(z) = 
  \begin{cases}
   F_+(z) & \im z > 0 \\ 
   F_-(z) & \im z < 0. 
  \end{cases}
\end{equation*}
We often use the representation 
\begin{align*}
 f(x) = \: & F(x+\mathrm{i}0) - F(x-\mathrm{i}0) \\ 
 = \: & F_+(x+\mathrm{i}0) - F_-(x-\mathrm{i}0). 
\end{align*}
If the limit 
\begin{equation*}
 \lim_{\epsilon\rightarrow 0+}
  \left\{ F(x+\mathrm{i}\epsilon) - F(x-\mathrm{i}\epsilon) \right\}
\end{equation*}
exists for $x\in I$, we define the value of the hyperfunction $f=[F]$ at
the point $x$ by this limit, and we do not define the value of $f$ if
the limit does not exist. 

We should remark that there is an ambiguity of the defining function
$F(z)$ of a hyperfunction $f(x)\in\mathscr{B}(I)$ up to a function
belonging to $\mathscr{O}(D)$. 
We mean that, if $F(z)\in\mathscr{O}(D\setminus I)$ defines a
hyperfunction $f(x)\in\mathscr{B}(I)$, a function $F(z)+\varphi(z)$ with 
$\varphi\in\mathscr{O}(D)$ also belongs to $\mathscr{O}(D\setminus I)$
and defines the same hyperfunction $f(x)\in\mathscr{B}(I)$, that is, 
\begin{equation*}
 f(x) = [F(z)] = [F(z) + \varphi(z)]. 
\end{equation*}

The addition $(f+g)(x)$ of two hyperfunctions $f(x)=[F(z)]$ and
$g(x)=[G(z)]$ are defined by 
\begin{equation*}
 (f+g)(x) := [F(z) + G(z)], 
\end{equation*}
and the multiplication  $cf(x)$ of a hyperfunction $f(x)=[F(z)]$ by a complex constant
$c$ is defined by
\begin{equation*}
 cf(x) := [cF(z)].
\end{equation*}
The space $\mathscr{B}(I)$ of hyperfunctions on an interval $I$ forms a
complex linear space by the definitions of the addition and the scalar
product. In addition, the multiplication $\varphi(x)f(x)$ of a
hyperfunction $f(x)=[F(z)]$ by a real analytic function $\varphi(x)$, 
that is, 
a real valued function $\varphi(x)$ on an interval $I$ which can be
extended to a complex function $\varphi(z)$ holomorphic on a complex
neighborhood $D$ of $I$, is defined by 
\begin{equation*}
 \varphi(x)f(x) := [\varphi(z)F(z)]. 
\end{equation*}
The above definitions are well-defined, that is, the definitions are not
dependent on the choice of the defining functions of hyperfunctions. 
For example, if $f(x) = [F(z)] = [\widetilde{F}(z)]$ and $g(x) = [G(z)]
= [\widetilde{G}(z)]$, we have 
\begin{equation*}
 f(x) + g(x) = [F(z) + G(z)] = [\widetilde{F}(z) + \widetilde{G}(z)]. 
\end{equation*}
The derivatives of a hyperfunction $f(x)=[F(z)]$ are defined by 
\begin{gather*}
 f^{\prime}(x) := [F^{\prime}(z)], \\ 
 f^{(n)}(x) := [F^{(n)}(z)] \quad ( \: n = 1, 2, \ldots \: ). 
\end{gather*}
Therefore, hyperfunctions are infinitely differentiable because
holomorphic functions are infinitely differentiable. 
The derivatives of a hyperfunction are also well-defined in the above
sense. 

\begin{figure}[h]
 \begin{center}
  \begin{tabular}{cc}
   \includegraphics[width=0.49\textwidth]{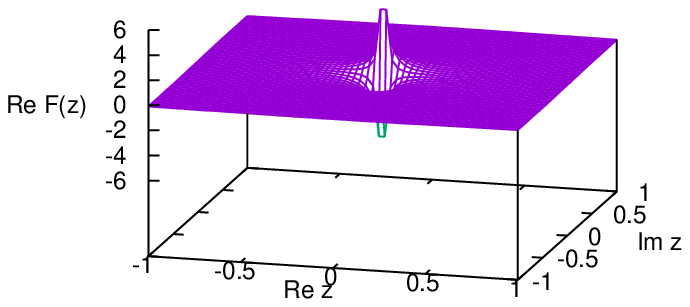} & 
       \includegraphics[width=0.49\textwidth]{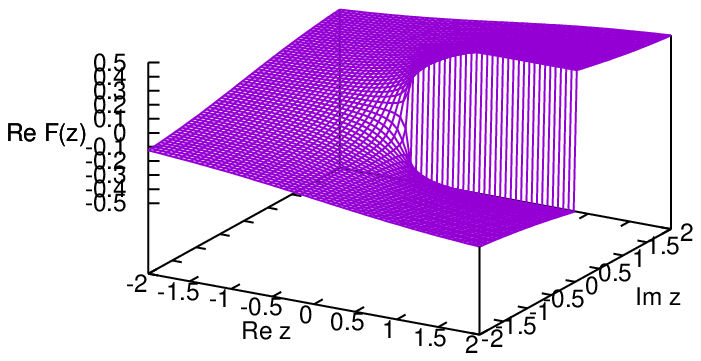}
       \\
   (a) $\delta(x)$ & (b) $Y(x)$
  \end{tabular}
 \end{center}
 \caption{The real part of the defining functions of (a) the delta
 function $\delta(x)$ and (b) the step function $Y(x)$. }
 \label{fig:hyperfunction-example}
\end{figure}

We here show some typical examples of hyperfunctions. 
The Dirac delta function $\delta(x)$ is given as a
hyperfunction by 
\begin{equation*}
 \delta(x) 
  = 
  \left[
    - \frac{1}{2\pi\mathrm{i}}\frac{1}{z}
  \right] 
  = 
  - \frac{1}{2\pi\mathrm{i}}
  \left( \frac{1}{x+\mathrm{i}0} - \frac{1}{x-\mathrm{i}0} \right), 
\end{equation*}
and the Heaviside step function 
\begin{equation*}
 Y(x) = 
  \begin{cases}
   1 & x > 0 \\ 
   0 & x < 0
  \end{cases}
\end{equation*}
is given as a hyperfunction by 
\begin{equation*}
 Y(x) = 
  \left[ -\frac{1}{2\pi\mathrm{i}}\log(-z) \right]
  = 
  - \frac{1}{2\pi\mathrm{i}}
  \left[ \log(-(x+\mathrm{i}0)) - \log(-(x-\mathrm{i}0))\right], 
\end{equation*}
where the complex logarithmic function $\log z$ is the branch such that 
$\log x$ is real valued if $x$ is a real positive number. 
Figure \ref{fig:hyperfunction-example} shows the graphs of the defining
functions of the delta function $\delta(x)$ and the step function
$Y(x)$. 
From these figures, we can understand visually that the difference
between the boundary values of the defining function gives a
hyperfunction. 
%
\section{Numerical Fourier transform}
\label{sec:numerical-Fourier-transform}
\begin{figure}[h]
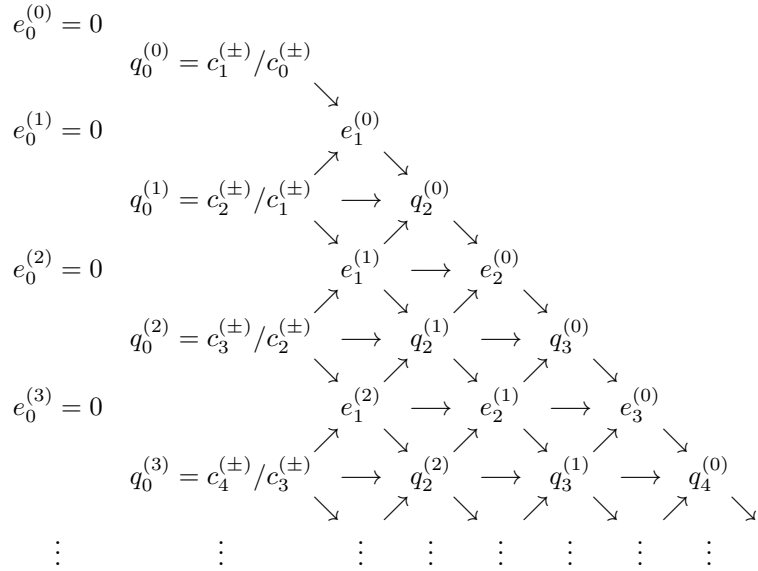

 \begin{equation*}
  \begin{array}{cc@{}c@{}c@{}c@{}c@{}c@{}c@{}c@{}c@{}c@{}c@{}c@{}c@{}c}
   e_0^{(0)} = 0 & \\
   & q_0^{(0)} = c_1^{(\pm)} / c_0^{(\pm)} \\ 
   & & \searrow \\ 
   e_0^{(1)} = 0 & & & e_1^{(0)} & \\ 
   & & \nearrow & & \searrow \\ 
   & q_0^{(1)} = c_2^{(\pm)} / c_1^{(\pm)} & & \longrightarrow & & q_2^{(0)} \\ 
   & & \searrow & & \nearrow & & \searrow \\ 
   e_0^{(2)} = 0 & & & e_1^{(1)} & & \longrightarrow & & e_2^{(0)}\\ 
   & & \nearrow & & \searrow & & \nearrow & & \searrow \\ 
   & q_0^{(2)} = c_3^{(\pm)} / c_2^{(\pm)} & & \longrightarrow & & q_2^{(1)} & & \longrightarrow & & q_3^{(0)} \\ 
   & & \searrow & & \nearrow & & \searrow & & \nearrow & & \searrow \\ 
   e_0^{(3)} = 0 & & & e_1^{(2)} & & \longrightarrow & & e_2^{(1)} & & \longrightarrow & &
    e_3^{(0)} \\ 
   & & \nearrow & & \searrow & & \nearrow & & \searrow & & \nearrow & &
    \searrow \\ 
   & q_0^{(3)} = c_4^{(\pm)} / c_3^{(\pm)} & & \longrightarrow & & q_2^{(2)} & & \longrightarrow & &
    q_3^{(1)} & & \longrightarrow & & q_4^{(0)} \\ 
   & & \searrow & & \nearrow & & \searrow & & \nearrow & & \searrow & &
    \nearrow & & \searrow
    \\ 
   \vdots & \vdots & & \vdots & & \vdots & & \vdots & & \vdots & &
    \vdots & & \vdots 
  \end{array}
 \end{equation*}
 \caption{The tableau of the quotient difference algorithm.}
 \label{fig:QD-algorithm}
\end{figure}

As mentioned in Section \ref{sec:introduction}, 
the Fourier transform of a function $f(x)$ is given as the hyperfunction 
\begin{equation}
 \label{eq:Fourier-transform3}
 \mathscr{F}[f](\xi) = [\mathfrak{F}_+(\zeta), \mathfrak{F}_-(\zeta)] = 
  \mathfrak{F}_+(\xi + \mathrm{i}0) - \mathfrak{F}_-(\xi + \mathrm{i}0)
\end{equation}
in hyperfunction theory, where 
\begin{gather}
 \label{eq:defining-function-0+}
 \mathfrak{F}_+(\zeta) = 
 \int_{-\infty}^{0}f(x)\mathrm{e}^{-2\pi\mathrm{i}\zeta x}\mathrm{d}x 
 \quad ( \: \im\zeta > 0 \: ) 
 \intertext{and}
 \label{eq:defining-function-0-}
 \mathfrak{F}_-(\zeta) = 
 - \int_0^{\infty}f(x)\mathrm{e}^{-2\pi\mathrm{i}\zeta x}\mathrm{d}x
 \quad ( \: \im\zeta < 0 \: ). 
\end{gather}
The function $\mathfrak{F}_+(\zeta)$ is holomorphic in the upper half plane 
$\im\zeta > 0$ 
and the function $\mathfrak{F}_-(\zeta)$ is holomorphic in the lower half plane 
$\im\zeta < 0$. 
Then, the Fourier transform $\mathscr{F}[f](\xi)$ is given as the
hyperfunction whose defining functions are $\mathfrak{F}_+(\zeta)$ and
$\mathfrak{F}_-(\zeta)$ by (\ref{eq:Fourier-transform3}). 
We remark that it is easy to compute the functions $\mathfrak{F}_{\pm}(\zeta)$ 
$( \: \pm\im\zeta > 0 \: )$ 
even if the integrand on the right-hand-side of
(\ref{eq:Fourier-transform}) 
is a slowly decaying oscillatory function  
because the integrands on the right-hand
sides of (\ref{eq:defining-function-0+}) and
(\ref{eq:defining-function-0-}) includes the exponentially decaying
factors $\exp(-2\pi|\im\zeta|x)$ 
Therefore we expect that we can compute the Fourier transform $\mathscr{F}[f](\xi)$ 
by the following way.
\begin{enumerate}
 \item We compute the defining functions $\mathfrak{F}_{\pm}(\zeta)$ 
       in $\mathbb{C}\setminus\mathbb{R}$. 
 \item We obtain the Fourier transform 
       \begin{math}
       \mathscr{F}[f](\xi) = 
	\mathfrak{F}_+(\xi+\mathrm{i}0) - \mathfrak{F}_-(\xi-\mathrm{i}0)
       \end{math} 
       by the analytic continuation of $\mathfrak{F}_{\pm}(\zeta)$
       onto the real axis $\mathbb{R}$.
\end{enumerate}
We call the method above the ^^ ^^ hyperfunction method''.
More exactly, the hyperfunction method for numerical Fourier transforms
is as follows. 
\paragraph{Computation of $\mathfrak{F}_{\pm}(\zeta)$}
We compute the defining functions $\mathfrak{F}_{\pm}(\zeta)$ in
$\mathbb{C}\setminus\mathbb{R}$ in Taylor series
\begin{gather}
 \label{eq:Taylor-series}
 \mathfrak{F}_{\pm}(\zeta) = 
 \sum_{n=0}^{\infty}c_n^{(\pm)}(\zeta - \zeta_0^{(\pm)})^n, 
 \intertext{%
 where $\zeta_0^{(\pm)}$ are given imaginary numbers such
 that $\pm\im\zeta_0^{(\pm)}>0$ and %
 }
 \begin{aligned}
  c_n^{(\pm)} = \: & 
  \frac{1}{n!}\mathfrak{F}_{\pm}^{(n)}(\zeta_0^{(\pm)}) 
  = 
  \pm\frac{1}{n!}\int_0^{\infty}(\pm 2\pi\mathrm{i}x)^n
  f(\mp x)\mathrm{e}^{\pm 2\pi\mathrm{i}\zeta_0^{(\pm)}x}
  \mathrm{d}x
  \\ 
  = \: & 
  \frac{1}{2\pi\eta n!}\int_0^{\infty}
  \left(\frac{\mathrm{i}x}{\eta}\right)^n
  f\left(-\frac{x}{2\pi\eta}\right)\mathrm{e}^{\mathrm{i}(\xi/\eta)x}
  \mathrm{e}^{-x}\mathrm{d}x
  \\ 
  & \hspace{40mm}( \: \zeta = \xi + \mathrm{i}\eta; \: n = 0, 1, 2, \ldots \: ).
 \end{aligned}
 \label{eq:Taylor-series-coefficient}
\end{gather}
It is easier to compute the coefficients $c_n$ than the oscillatory
integral on the right-hand side of (\ref{eq:Fourier-transform}) 
because the integrands on the right-hand side of
(\ref{eq:Taylor-series-coefficient}) involves the exponentially
decaying factor $\mathrm{e}^{-2\pi|\im\zeta_0^{\pm}|x}$. 
We can compute the integrals appearing in
(\ref{eq:Taylor-series-coefficient}) easily using conventional
quadrature formula, for example, the DE rule.
\paragraph{Analytic continuation of $\mathfrak{F}_{\pm}(\zeta)$}
Next, we obtain the analytic continuation of $\mathfrak{F}_{\pm}(\zeta)$ 
onto the real axis $\mathbb{R}.$
For this purpose, we transform the Taylor series
(\ref{eq:Taylor-series}) 
into the continued fractions 
\begin{equation*}
 \mathfrak{F}_{\pm}(\zeta) = 
  \cfrac{a_1^{(\pm)}}{1 + \cfrac{a_2^{(\pm)}(\zeta-\zeta_0^{(\pm)})}{1 +
  \cfrac{a_3^{(\pm)}(\zeta-\zeta_0^{(\pm)})}{1 + \ddots}}},
\end{equation*}
We expect that we can get the analytic continuation of
$\mathfrak{F}_{\pm}$ by the continued fraction 
since, in general, the convergence regions of the continued fractions are wider than
the convergence disk of the Taylor series (\ref{eq:Taylor-series}) \cite{Henrici1977}. 
The coefficients $a_n^{(\pm)}$ of the continued fractions are obtained
by the quotient-difference algorithm \cite{Henrici1977} as follows. 
We compute the series 
$\{e_k^{(n)}|\:n=0,1,2,\ldots;k=0,1,2,\ldots\}$ and 
$\{q_k^{(n)}|\:n=0,1,2,\ldots;k=1,2,\ldots\}$ by 
\begin{gather*}
 e_0^{(n)} = 0, \quad  
 q_1^{(n)} = \frac{c_{n+1}^{(\pm)}}{c_n^{(\pm)}} \quad ( \: n = 0, 1, 2, \ldots \: ), 
 \\ 
 \begin{split}
  e_k^{(n)} = 
  q_k^{(n+1)} - q_k^{(n)} + e_{k-1}^{(n+1)}, \quad 
  q_{k+1}^{(n+1)} = 
  \frac{e_k^{(n+1)}}{e_k^{(n)}}q_k^{(n+1)} 
  \\
  ( \: n = 0, 1, 2, \ldots; \: k = 1, 2, \ldots \: ),
 \end{split}
\end{gather*}
and then we obtain the coefficients $a_n$ by
\begin{equation*}
 a_1^{(\pm)} = c_0^{(\pm)}, \quad 
  a_{2k}^{(\pm)} = - q_k^{(0)}, \quad 
  a_{2k+1}^{(\pm)} = - e_k^{(0)} \quad ( \: k = 1, 2, \ldots \: ).
\end{equation*}
The numbers $e_k^{(n)}$ and $q_k^{(n)}$ are generated as shown in the
tableau of Figure \ref{fig:QD-algorithm}. 
It is known that the quotient-difference algorithm is numerically
unstable. 
Therefore we carry out the computation of the algorithm using multiple precision
arithmetics. 
If some of the coefficients $c_k^{(\pm)}$ are zeros, we cannot carry out
the algorithm. 
In that case, we change $\zeta_0^{(\pm)}$ so that all the coefficients
$c_k^{(\pm)}$ are not zeros. 
%
\section{Numerical examples}
\label{sec:example}
\begin{figure}[h]
 \begin{center}
  \begin{tabular}[t]{cc}
   \includegraphics[width=0.49\textwidth]{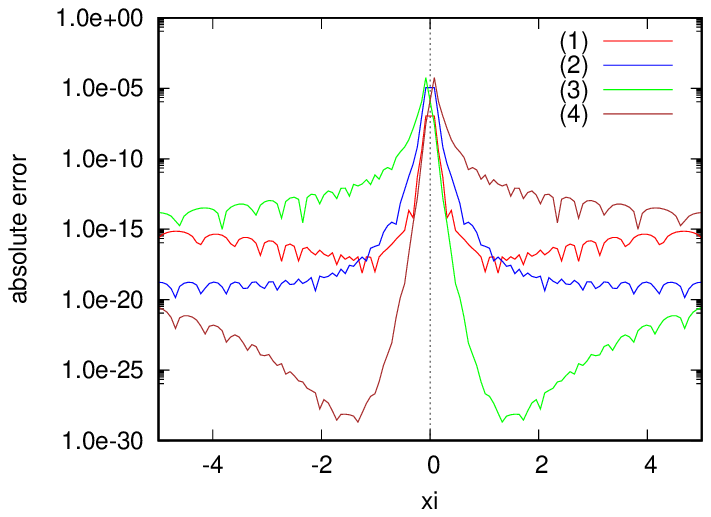} & 
       \includegraphics[width=0.49\textwidth]{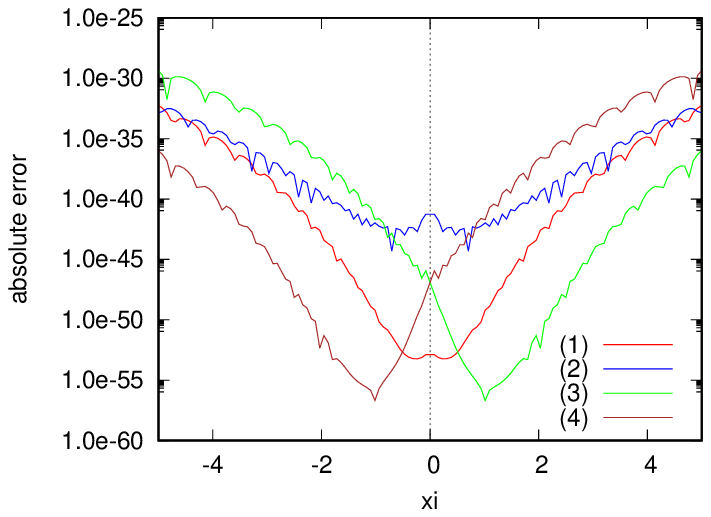}
       \\
   (i) $f(x) = (1 + x^2)^{-1}$ & (ii) $f(x) = \tanh(\pi x)$
       \\ 
   \includegraphics[width=0.49\textwidth]{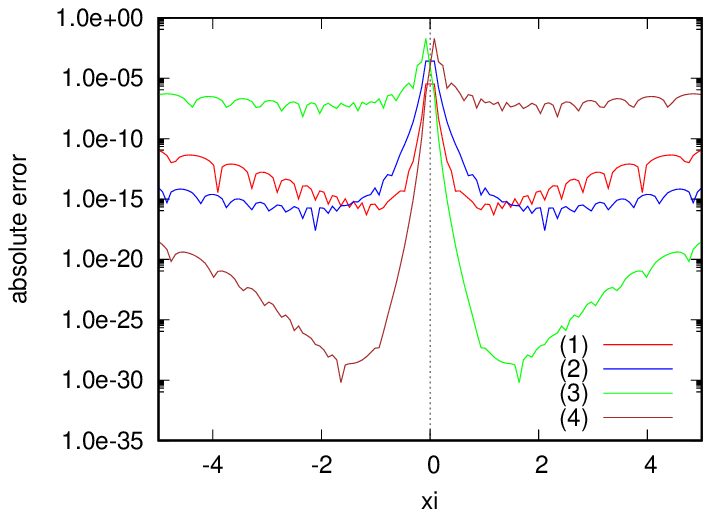} & 
       \includegraphics[width=0.49\textwidth]{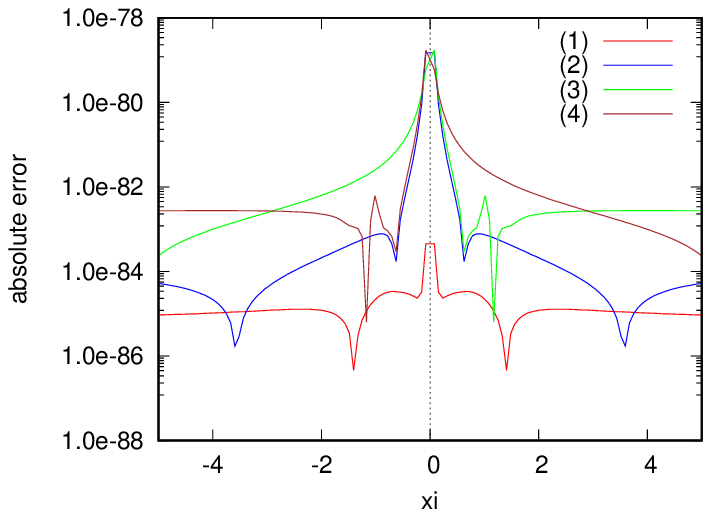}
       \\ 
   (iii) $f(x) = \log|x|$ & (iv) $f(x) = |x|$
  \end{tabular}
 \end{center}
 \caption{The absolute errors of the numerical Fourier transforms of the
 functions (\ref{eq:function-example1}) by the hyperfunction method with 
 the center of the Taylor series in (\ref{eq:Taylor-series}) taken as
 (1) $\zeta_0^{(\pm)} = \pm\mathrm{i}$, (2) $\zeta_0^{(\pm)} = \pm 2\mathrm{i}$,
 (3) $\zeta_0^{(\pm)} = 1\pm\mathrm{i}$ and $\zeta_0^{(\pm)} = - 1 \pm
 \mathrm{i}$.}
 \label{fig:example1}
\end{figure}

We here show some numerical examples which show the effectiveness of the
proposed method. 
All the computations were carried out by using programs coded in C++
with 100 decimal digit precision working by the multiple precision
arithmetic library {\it exflib} \cite{Fujiwara-exflib}. 
We computed the Fourier series $\mathscr{F}[f](\xi)$ for 
\begin{equation}
 \label{eq:function-example1}
 f(x) = 
  \begin{cases}
   \mathrm{(i)} \: & 1 / (1 + x^2) \\ 
   \mathrm{(ii)} \: & \tanh(\pi x) \\ 
   \mathrm{(iii)} \: &  \log|x| \\ 
   \mathrm{(iv)} \: & |x|
  \end{cases}
\end{equation}
by the hyperfunction method. 
The exact expression of the Fourier transforms for the functions
(\ref{eq:function-example1}) are known as follows.
\begin{align*}
 \mathrm{(i)} \quad & \mathscr{F}[(1+x^2)^{-1}](\xi) =
 \pi\mathrm{e}^{-2\pi|\xi|}, \\ 
 \mathrm{(ii)} \quad & \mathscr{F}[\tanh(\pi x)](\xi) = - \mathrm{i}\cosech(\pi\xi), \\ 
 \mathrm{(iii)} \quad & \mathscr{F}[\log|x|](\xi) = - \gamma\delta(\xi) - \frac{1}{2|\xi|},
 \\ 
 \mathrm{(iv)} \quad & \mathscr{F}[|x|](\xi) = - \frac{1}{2(\pi\xi)^2}.
\end{align*}
Figure \ref{fig:example1} shows the absolute errors of the numerical Fourier
transforms by the hyperfunction method with the center of Taylor series 
\begin{math}
 \zeta_0^{(\pm)} = \pm\mathrm{i}, \ \pm 2\mathrm{i}, \ 
 1\pm\mathrm{i}
\end{math} 
and $-1\pm\mathrm{i}$, 
and Table \ref{tab:example1} shows the numbers of the evaluations of
$f(x)$ in (\ref{eq:function-example1}) in computing the Fourier
transforms of the functions $f(x)$ by the hyperfunction method. 
From the figures, the hyperfunction method works well, especially, for 
$f(x) = \tanh(\pi x)$ and $|x|$. 
\begin{table}[h]
 \caption{The numbers of the evaluations of $f(x)$ in computing the
 Fourier transform of the functions (\ref{eq:function-example1}) by the
 hyperfunction method.}
 \begin{center}
  \begin{tabular}[t]{ccccc}
   \hline \rule{0pt}{12pt}
   $\zeta_0^{(\pm)}$ & $\pm\mathrm{i}$ & $\pm 2\mathrm{i}$ &
	       $1\pm\mathrm{i}$ & $-1\pm\mathrm{i}$ 
		   \\
   \hline \rule{0pt}{12pt}
   (i) & 1420 & 720 & 2820 & 2820 \\ 
   (ii) & 1330 & 666 & 2838 & 2838 \\ 
   (iii) & 1430 & 714 & 2838 & 2838 \\ 
   (iv) & 1332 & 668 & 2646 & 2646 \\ 
   \hline
  \end{tabular}
 \end{center}
 \label{tab:example1}
\end{table}

We remark the interesting fact that 
we can obtain Fourier transforms by the hyperfunction method without
computing oscillatory integrals. In fact, if we take $\zeta_0^{(\pm)}$
as a purely imaginary number $\xi = \mathrm{i}\eta$, 
the integrals in (\ref{eq:Taylor-series-coefficient}) which gives the
coefficients of the Taylor series of $\mathfrak{F}_{\pm}(\zeta)$ become 
\begin{equation*}
 c_n^{(\pm)} = 
  \frac{1}{2\pi\eta n!}\int_0^{\infty}
  \left(\frac{\mathrm{i}x}{\eta}\right)^n
  f\left(-\frac{x}{2\pi\eta}\right)\mathrm{e}^{-x}\mathrm{d}x, 
\end{equation*}
which include no oscillatory function. 

We compared the hyperfunction method with Sugihara's method and the DE-type
formula by Ooura and Mori.
We computed the Fourier transform $\mathscr{F}[f](\xi)$ for $f(x)$ given
in (\ref{eq:function-example1}) with $\xi = 1$ 
by our method and the two previous methods. 
Table \ref{tab:example2} shows the numbers of the evaluations of $f(x)$ 
used for numerical integration and the errors of the methods. 
From the table, the our method is superior to Sugihara's method, 
and it is competitive with Ooura and Mori's method in some examples. 
Besides, we remark that our method gives Fourier transform as a function
while the two previous methods give Fourier transform as an integral or
a number. 
We means that, in our method,  we can use the same coefficients $a_n$ 
of the continued fraction for Fourier
transforms $\mathscr{F}[f](\xi)$ with many $\xi$'s once the coefficients
$a_n$ are obtained. 
On the other hand, in the two previous methods, we have to compute
numerical integrals again every time the value $\xi$ is changed. 
\begin{table}[h]
 \caption{The numbers $N$ of the evaluations of $f(x)$ and the errors of the
 hyperfunction method, Sugihara's method and Ooura-Mori's method applied
 to the Fourier transforms $\mathscr{F}[f](\xi=1)$ of the functions $f(x)$ in
 (\ref{eq:function-example1}).}
 \begin{center}
  \begin{tabular}[t]{cccc}
   \hline \rule{0pt}{12pt}
   $f(x)$ & method & $N$ & error \\
   \hline \rule{0pt}{12pt} 
   & hyperfunction $(\:\zeta_0^{(\pm)}=\pm\mathrm{i}\:)$ & 1420 &
	       $8.1\times 10^{-18}$ 
	       \\
   & hyperfunction $(\:\zeta_0^{(\pm)}=\pm 2\mathrm{i}\:)$ & 710 & 
	       $1.0\times 10^{-16}$
	   \\
   (i) & hyperfunction $(\:\zeta_0^{(\pm)}=1\pm\mathrm{i}\:)$ & 2820 &
	       $3.1\times 10^{-27}$
   \\ 
   & Sugihara & 17808 & $8.1\times 10^{-22}$
       \\ 
   & Ooura-Mori & 1010 & $2.2\times 10^{-28}$
	       \\
   \hline \rule{0pt}{12pt}
   & hyperfunction $(\:\zeta_0^{(\pm)}=\pm\mathrm{i}\:)$ & 1330 &
	       $4.9\times 10^{-50}$ 
	       \\ 
   & hyperfunction $(\:\zeta_0^{(\pm)}=\pm 2\mathrm{i}\: )$ & 666 &
	       $7.4\times 10^{-43}$
	       \\ 
   (ii) & hyperfunction $(\:\zeta_0^{(\pm)}=1\pm\mathrm{i}\:)$ & 2642 &
	       $2.1\times 10^{-56}$
	       \\
   & Sugihara & 17156 & $7.9\times 10^{-21}$
       \\ 
   & Ooura-Mori & 1892 & $1.5\times 10^{-46}$ 
	       \\ 
   \hline \rule{0pt}{12pt}
   & hyperfunction $(\:\zeta_0^{(\pm)}=\pm\mathrm{i}\:)$ & 1430 &
	       $4.2\times 10^{-16}$ 
	       \\ 
   & hyperfunction $(\:\zeta_0^{(\pm)}=\pm 2\mathrm{i}\:)$ & 714 &
	       $9.8\times 10^{-28}$ 
	       \\ 
   (iii) & hyperfunction $(\:\zeta_0^{(\pm)}=1\pm\mathrm{i}\:)$ & 2838 &
	   $2.8\times 10^{-28}$ 
	   \\ 
   & Sugihara & 17916 & $2.0\times 10^{-20}$
       \\ 
   & Ooura-Mori & 1020 & $6.3\times 10^{-32}$ 
	   \\ 
   \hline
   & hyperfunction $(\:\zeta_0^{(\pm)}=\pm\mathrm{i}\:)$ & 1332 &
	       $2.2\times 10^{-85}$
	   \\ 
   & hyperfunction $(\:\zeta_0^{(\pm)}=\pm 2\mathrm{i}\:)$ & 668 &
	       $7.5\times 10^{-84}$
	   \\ 
   (iv) & hyperfunction $(\:\zeta_0^{(\pm)}=1\pm\mathrm{i} \: )$ & 2646
	   & $6.3\times 10^{-83}$
       \\ 
   & Sugihara & 17182 & $3.8\times 10^{-20}$
	   \\ 
   & Ooura-Mori & 15198 & $1.5\times 10^{-80}$
   \\
   \hline
  \end{tabular}
 \end{center}
 \label{tab:example2}
\end{table}
\section{Concluding remarks}
\label{sec:concluding-remarks}
In this paper, we proposed a numerical method for obtaining Fourier
transforms based on hyperfunction theory. 
In hyperfunction theory, a Fourier transform is given as a
hyperfunction, the difference of the real axis of analytic functions
which are called the defining functions of the hyperfunction, 
and the defining functions of a Fourier transform are given by integrals
including the integrand of the desired Fourier transform and 
exponentially decaying factors. 
In our method, we compute the defining functions in the upper or lower
half complex plane and obtain the Fourier transform by the analytic
continuations of the defining functions onto the real axis. 
Numerical examples show that our method is effective and competitive
with the previous methods. 

Problems for future study are as follows. 
\begin{itemize}
 \item How should we choose the centers $\zeta_0^{(\pm)}$ of the
       Taylor series of the defining functions
       $\mathfrak{F}_{\pm}(\zeta)$?  
       As shown in the numerical examples, the accuracy of our method
       depends on the choice of $\zeta_0^{(\pm)}$, and it is a crucial
       problem which points $\zeta_0^{(\pm)}$ are the best. 
 \item Theoretical error estimate of our method is an important
       problem. 
 \item We employed an analytic continuation by transforming the defining
       functions into continued fractions. However, this process is
       expensive because we use the quotient-difference algorithm here
       and we use multiple precision arithmetics due to the numerical
       instability of the quotient-difference algorithm. 
       Therefore, we need a numerically stable method of analytic
       continuation or transformation of an analytic function into a
       continued fraction. 
\end{itemize}
\section*{Acknowledgements}
This work is supported by JSPS KAKENHI Grant Number JP16K05267. 

%
%

\begin{thebibliography}{00}
\bibitem{Fujiwara-exflib} H. Fujiwara, 
	Exflib information, \\
	http://www-an.acs.i.kyoto-u.ac.jp/\~{}fujiwara/exflib/.
\bibitem{Graf2010} U. Graf, 
	Introduction to Hyperfunctions and Their Integral Transforms 
	--- An Applied and Computational Approach, Birkh\"{a}user,
	Basel, 2010. 
\bibitem{Henrici1977} P. Henrici, 
	Applied and Computational Complex Analysis, Vol. 2, 
	John Wiley \& Sons, New York, 1977. 
\bibitem{Kaneko1988} A, Kaneko, 
	Introduction to Hyperfunctions,  
	Kluwer Academic Publications, Boston, 1988. 
\bibitem{Mori1972} M. Mori, 
	Numerical analysis and hyperfunction theory, 
	^^ ^^ Kokyuroku'', Res. Inst. Math. Sci. Kyoto Univ. 
	145 (1972) 1--11 (in Japanese).
\bibitem{OgataHirayama2018} H. Ogata and H. Hirayama, 
	Numerical integration based on hyperfunction theory, 
	J. Comput. Appl. Math. 327 (2018) 243--259.
\bibitem{OouraMori1991} T. Ooura and M. Mori, 
	The double exponential formula for oscillatory functions over
	the half infinite interval, 
	J. Comput. Appl. Math. 38 (1991) 353--360.
\bibitem{Sato1959} M. Sato, 
	Theory of hyperfunctions, J. Fac. Sci. Univ. Tokyo, Sect. 1A
	Math. 8 (1959) 139--193.
\bibitem{Sugihara1987} M. Sugihara, 
	Methods of numerical integration of oscillatory functions by the
	DE-formula with the Richardson extrapolation, 
	J. Comput. Appl. Math. 17 (1987) 47--68.
\bibitem{TakahasiMori1974} H. Takahasi and M. Mori, 
	Double exponential formulas for numerical integration, 
	Publ. RIMS, Kyoto Univ. 339 (1978) 721--741.
\bibitem{TodaOno1978} H. Toda and H. Ono, 
	Some remarks for efficient usage of the double exponential
	formulas (in Japanese), Kokyuroku, RIMS, Kyoto Univ. 
	339 (1978) 74--109. 
\end{thebibliography}
\end{document}